\documentclass[12pt]{article}
\usepackage{amsfonts}
\begin{document}
\title{\bf  Existence of minimizers of functionals involving the
fractional gradient in the abscence of compactness, symmetry and
monotonicity
\author{ H. Hajaiej\\
  }}
\date{}
\maketitle
\begin{abstract}
We establish general assumptions under which a constrained
variational problem involving the fractional gradient and a local
nonlinearity admits minimizers.
\end{abstract}
\section{Intoduction}

For a prescribed number $c > 0$ and $0 < s < 1$,
we consider the
following constrained minimization problem :
$$\inf \{J(u) : u \in S_c\} = I_c\eqno{(1.1)}$$
\begin{eqnarray*}
J(u) &=& \frac{1}{2} \int\; |\nabla_s u|^2 - \int F(x,u),\\
|\nabla_s u|^2_2 &=& \int |\nabla_s u|^2 = C_{N,s}
\int\int\;\frac{|u(x) - u(y)|^2}{|x-y|^{N+2s}}dx dy,
\end{eqnarray*}
$F$ is a carath\'eodory function, and
$$S_c = \{u \in H^s(\mathbb{R}^N) : \int u^2 = c^2\}.$$
Under some additional regularity assumptions on $F$, solutions of
(1.1) satisfy the following fractional elliptic equation :
$$\Delta^s u + f(x,u) + \lambda u = 0\eqno{(1.2)}$$
where $F(x,t) = \displaystyle{\int^t_0}f(x,p) dp$ and $\lambda$ is a
Lagrange multiplier. Solutions of (1.1) can also be viewed as
standing waves of the following nonlinear fractional Schr\"odinger
equation
$$\left\{ \begin{array}{l}
i\partial_t \Phi(t,x) + f(x, |\Phi|) + \Delta^s_{xx}\Phi = 0\\
\Phi(0,x) = \Phi^0(x).
\end{array}\right.\eqno{(1.3)}$$
Despite the importance of (1.2) and (1.3) in many domains, there are
only results the particular cases : $N = 1,  s = \frac{1}{2}$ and
$f(x,s) = s^\alpha , [1,2]$. Let us point out that when $N = 3, s =
\frac{1}{4}$, (1.3) models water waves,  semilunar heart valve
vibrations and neural systems. When $s = \frac{3}{4}$, it governs
water waves with surface tension, [5]. More generally, equations
(1.2) and (1.3) arise in numerous models from mathematical physics,
mathematical biology, finance, inhomogenous porous material,
geology, hydrology, dynamics of earthquakes, bioegineering, chemical
engineering, neural networks and medicine,
[5,6] and references therein.\\
In this paper, we address the question of existence of minimizers of
(1.1) in the absence of compactness,  symmetry and monotonicity .
This considerably extends the main result obtained by the author in
[5], where the integrand $F$ has a nice combination of monotonicity
and symmetry properties, which enabled us to obtain the compactness
of Schwarz minimizing sequences. In the present work, we will prove
the above property for any minimizing sequence of (1.1) without
requiring any symmetry or monotonicity properties of the
integrand.\\
Our main result is :\vspace{3mm}\\
{\bf Theorem 1.1.} Suppose that the function $F : \mathbb{R}^N
\times \mathbb{R}  \rightarrow \mathbb{R}$ is a Carath\'eodory
function verifying :

(F0) $\forall\; x \in \mathbb{R}^N, t \in \mathbb{R}, \; \exists\;
A, A' > 0$ and $0 < \ell < \frac{4s}{N}$ such that :
$$0 \leq F(x,t) \leq A(t^2 + |t|^{\ell+2})$$
and
$$0 \leq \partial_2 F(x,t) \leq A'(|t| + |t|^{\ell+1})$$

(F1) $\exists\; \Delta > 0, S > 0, R > 0, \alpha > 0$ $p\in [0,2)$
such that :
$$F(x,t) > \Delta |x|^{-p} |t|^\alpha \mbox{ for } |x|\geq R, |t| < S,$$
where
$$N+2s > \frac{N}{2} \alpha+p,$$

(F2) $F(x,\theta t) \geq \theta^2 F(x,t) \; \forall\; x \in
\mathbb{R}^N, t \in \mathbb{R}$ $\theta \geq 1$.\\ There exists a
periodic function $F^\infty (x,t)$ (i.e $\exists\; z \in
\mathbb{Z}^N$ such that $F^\infty(x+z,t) = F^\infty(x,t)$ $\forall\;
x \in \mathbb{R}^N, t \in \mathbb{R})$ satisfying (F1) such that :

(F3) There exists $0 < \beta < \frac{4s}{N}$ such that
$\displaystyle{\lim_{|x| \rightarrow \infty}\;
\frac{F(x,t)-F^\infty(x,t)}{t^2 + |t|^{\beta+2}}} =0$ uniformly for
any $t$.

(F4) There exists $B , B'$ and $0 < \gamma < \ell < \frac{4s}{N}$
such that
$$0 \leq F^\infty (x,t) \leq B(|t|^{\gamma +2}+ |t|^{\ell+2})$$
and
$$0 \leq \partial_2 F^\infty (x,t) \leq B'(|t|^{\gamma+1} + |t|^{\ell+1})$$
$\forall\; x \in \mathbb{R}^N, t \in \mathbb{R}$.

(F5) There exists $\sigma \in (0, \frac{4s}{N})$ such that
$$F^\infty (x, \theta t) \geq \theta^{\sigma+2} F^\infty(x,t)$$
$\forall\; \theta\geq 1$, $x \in \mathbb{R}^N$ and $t \in
\mathbb{R}$.

(F6) $F^\infty(x,t) \leq F(x,t)\quad \forall\; x \in \mathbb{R}^N, t
\in \mathbb{R}$,  with strict inequality in a measurable set having
a positive Lebesgue measure.

Then there exists $u_c \in S_c$ such that
$$J(u_c) = I_c.$$
{\bf Theorem 1.2} If (F1) holds true for $F^\infty$, (F4) and (F5)
are satisfied, then there exists $u_c \in S_c$ such that
$$J^\infty(u_c) = I^\infty_c, \mbox{ where } J^\infty(u) = \frac{1}{2}
\int |\nabla_s u|^2 - \int F^\infty(x,u)$$ and
$$I^\infty_c = \inf \{J^\infty (u) : u \in S_c\}.\eqno{(1.4)}$$
Our proofs of the above results are based on a variant of the
breakthrough concentration-compactness principle (appendix).\\
Our line of attack consists of the following steps :

In order to prove that vanishing cannot occur, it is sufficient to
show the strict negativity of the value of the infinimum (Lemma
3.2). Then, to rule out dichotomy, we will first prove that the
minimization problem (1.4) is achieved \hfill (S1)\\
 and that :
$$I_c < I^\infty_c\quad \forall\; c > 0\eqno{(S2)}$$
$$I_c \leq I_{c-a} + I_a\quad \forall\; a \in (0,c)\eqno{(S3)}$$
(S2) and (S3) imply the strict subadditivity inequality
$$I_c < I^\infty_{c-a} + I_a \quad \forall\; a \in (0,c)\eqno{(S4)}$$
On the other hand, we will prove that thanks to our assumptions on
$F$, we certainly have for any minimizing sequence $(u_n)$ of (1.1)
that :
$$J(u_n) \geq J(u_{n,1}) + J^\infty(u_{n,2})-g(\delta)\eqno{(S5)}$$
where $g(\delta)\rightarrow 0$ as $\delta \rightarrow 0$.\\
The latter requires of course a deep and subtle study of the
functionals $J$ and $J^\infty$ (Lemma 3.1).\\
Finally the continuity of $I_c$ and $I^\infty_c$ enables us to
deduce that (S5) implies the following inequality :
$$I_c \geq I_a + I_{c-a}^\infty.\eqno{(S6)}$$
(S4) together with (S6) yield to a contradiction.\\
Once one knows that compactness is the only plausible alternative,
the strict inequality (S2) will be very helpful to conclude that any
minimizing sequence of (1.1) is compact (up to a subsequence). These
issues were heuristically discussed in the classical setting in the
seminal paper of Lions [7].
\section{Notations}

$\bullet\quad$ $N \in \mathbb{N}^\ast, 0 < s < 1$ and $N \geq 2s$.\\
$\bullet\quad$ $A$ constant $C$ can vary from line to line, we will
keep the same notation for it.\\
$\bullet\quad$ The norm of $L^p(\mathbb{R}^N)$ is denoted by
$|\;|_p$ or $|\;|_{L^p}$\\
 $\bullet\quad$ $H^s(\mathbb{R}^N) = \{u
\in L^2(\mathbb{R}^N)\; \displaystyle{\int}(1+|\xi|^{2s})|
\mathcal{F}u(\xi)|^2 d\xi < \infty\}$ where $\mathcal{F}$ denotes
the Fourier transform,  which is equivalent to
$$H^s(\mathbb{R}^N) = H^s = \{u \in L^2(\mathbb{R}^N) :
\frac{|u(x)-u(y)|}{|x-y|^{N/2+s}}
\in L^2 (\mathbb{R}^N \times \mathbb{R}^N)\}$$ endowed with the
natural norm :
$$|u|_{H^s} = \left(\int |u|^2 + \int\int \frac{|u(x)-u(y)|^2}{|x-y|^{N+2s}}dxdy
\right)^{1/2}$$
$$|\nabla_s u|^2_2 = C_{N,s} \int\int\frac{|u(x)-u(y)|^2}{|x-y|^{N+2s}}dxdy.$$
$2^\ast_s = \displaystyle{\frac{2N}{N-2s}}$ if $N > 2s$ and
$2^\ast_s = \infty$ if $N = 2s$.\\
$H^{-s}(\mathbb{R}^N) = H^{-s}$ is the dual space of $H^s$.

In an integral where no domain of integration is indicated, $t$ is
to be understood that the integral extends over the whole space
\section{Proof of the main result}
{\bf Lemma 3.1} If $F$ satisfies (F0) , then
\begin{itemize}
\item[(i)] \begin{itemize}
\item[a)] $J \in C^1(H^s,\mathbb{R})$ and there exists a constant $D >
0$ such that :
$$|J'(u)|_{H^{-s}} \leq D(|u|_{H^s} +|u|_{H^s}^{1 + \frac{4s}{N}})$$
for any $u \in H^s$.
\item[b)] $J^\infty \in C^1(H^s, \mathbb{R})$ and there exists a
constant $D' > 0$ such that :
\end{itemize}
$$|J^{'\infty}(u)|_{H^{-s}} \leq D'(|u|_{H^s}
+ |u|_{H^s}^{1+ \frac{4s}{N}})$$ for any $u \in H^s$.
\item[(ii)] $J(u) \geq A_1 |\nabla_s u|^2_2 - A_2c^2 - A_3
c^{(1-\sigma)(\ell+2)q}$\\
$J^\infty(u) \geq B_1 |\nabla_s u|^2_2 - B_2
c^{(1-\sigma_1)(\beta+2)q_1} B_3 c^{(1-\sigma)(\ell+2)q}\; \quad
\forall\; u \in S_c .$\\ $(\sigma, \sigma_1, q$ and $q_1$ will be
given below).
\item[(iii)] \begin{itemize}
\item[a)] $I_c > - \infty$ and any  sequence of (1.1) is bounded in
$H^s$.
\item[b)] $I^\infty_c > - \infty$ and any minimizing sequence of
(1.4) is bounded in $H^s$.
\end{itemize}
\item[(iv)] $c \mapsto I_c$ and $c \mapsto I^\infty_c$ are
continuous on $(0, \infty)$.
\end{itemize}
{\bf Proof : } Let $\varphi : \mathbb{R} \rightarrow \mathbb{R}$ be
the function defined by :
$$\left\{\begin{array}{ll}
\varphi(t) = 1 &\mbox{ if } |t| < 1\\
\varphi(t) = -|t|+2 &\mbox{ if } 1 \leq |t| \leq 2\\
\varphi(t) = 0 &\mbox{ if } |t| > 2
\end{array}\right.$$
$$\partial^1_2 F(x,t) = \varphi(t) \partial_2 F(x,t)\quad \mbox{ and }$$
$$|\partial^1_2 F(x,t)| \leq A(1+2^{\ell+1})|t|\eqno{(3.1)}$$
$$\partial^2_2 F(x,t) = (1-\varphi(t))\partial_2 F(x,t)$$
$$|\partial^2_2 F(x,t)| \leq 2A |t|^{1 + \frac{4s}{N}}\eqno{(3.2)}$$
Let $$p = \left\{ \begin{array}{ll} \frac{2N}{N+2s} &\mbox{ for } N
> 2s\\
\frac{4}{3} &\mbox{ if } N = 2s
 \end{array}\right.$$
 and $q = (1+ \frac{4s}{N})p$.\\
(3.1) and (3.2) imply that
 $\partial^1_2 F(x,.) \in C(L^2,L^2)$ and $\partial^2_2 F(x,.) \in C
 (L^q,L^p)$ and there exists a constant $K > 0$ such that :
 $$|\partial^1_2 F(x,u)|_2 \leq K |u|_2\quad \forall\; u \in L^2$$
 $$|\partial^2_2 F(x,u)|_p \leq K |u|_q^{1 + \frac{4s}{N}}\quad
 \forall\; u \in L^q.$$
 Noticing that $H^s$ is continuously embedded in $L^q$ since
 $q \in [2, \frac{2N}{N-2s}]$ for $N > 2s$ and $q \in [2,\infty)$ for
 $N = 2s$ and $L^p$ is continuously embedded in $H^{-s}$ since $p' \in [2,
 \frac{2N}{N-2s}]$ for  $N > 2s$ and $p' \in [2, \infty)$ for $N =
 2s$. We can assert that :
 $$\partial_2^1 F(x,.) + \partial^2_2 F(x,.) \in C (H^s_, H^{-s}), $$
 and there exists a constant $C > 0$ such that
 $$|\partial_2 F(x,u)|_{H^{-s}} \leq C \{|u|_{H^s} +
  |u|_{H^s}^{1 + \frac{4s}{N}}\}\eqno{(3.3)}$$
 for all $u \in H^s$.\\
 On the other hand :
 $$\int F(x,u) \leq A(|u|^2_2 + |u|^{\ell+2}_{\ell+2}) \leq C(|u|^2_{H^s}
  + |u|^{\ell+2}_{H^s})$$
 , which implies that $J \in C^1(H^s,\mathbb{R})$ by standard
 arguments of differential calculus. \\
 Therefore,
 $$|J'(u)|_{H^{-s}} \leq C\{|u|_{H^s} + |u|_{H^s}^{1 + \frac{4s}{N}}\}\quad
 \forall\; u \in H^s.$$

\begin{itemize}
\item[(i)] b) can be easily deduced following the same steps which yield to
similar estimates as (3.1) and (3.2)

\item[(ii)] These estimates were obtained in [4].

\item[(iii)] is a direct consequence of (ii)\end{itemize}
{\bf Proof of (iv)}

Consider $c > 0$ and a sequence $\{c_n\} \subset (0, \infty)$ such
that $c_n \rightarrow c$. For any $n \in \mathbb{N}$, there exists
$u_n \in S_{c_n}$ such that $I_{c_n} \leq J(u_n) \leq I_{c_n} +
\frac{1}{n}$.\\
By (iii), there exists $K > 0$ such that $|u_n|_{H^s} \leq K$  for
all $n \in \mathbb{N}$. Setting $w_n = \frac{c}{c_n} u_n$, we have
that $w_n \in S_c$ and $|u_n-w_n|_{H^s} \leq |1-\frac{c}{c_n}|\;
|u_n|_{H^s} \leq K|1- \frac{c}{c_n}|$ for any $n \in \mathbb{N}$.\\
Therefore, there exists $n_1$ such that $|u_n-w_n|_{H^s} \leq K+1$
for all $n \geq n_1$. By part (i), there exists a constant $L(K) >
0$ such that $\|J'(u)\|_{H^{-s}} \leq L(K)$ for all $u \in H^s$ such
that $|u|_{H^s} \leq 2K+1$.\\
So for all $n \geq n_1$ :
\begin{eqnarray*}
|J(w_n)-J(u_n)| &=& |\int^1_0\; \frac{d}{dt} J(t w_n +
(1-t)u_n)dt|\\
&\leq& \sup_{\|u\|_{H^s} \leq 2 K+1}
\|J'(u)\|_{H^{-s}}\|u_n-w_n\|_{H^s}\\
&\leq& L(K)K|1 - \frac{c}{c_n}|
\end{eqnarray*}
and so  liminf $I_{c_n} \geq I_c$.\hfill (3.4)

On the other hand there exists a sequence $\{u_n\} \subset S_c$ such
that $J(u_n) \rightarrow I_c$ and  thus by (iii), we can find $K >
0$ such that $|u_n|_{H^s} \leq K$. $w_n = \frac{c_n}{c} u_n$. As
above, we can write $w_n \in S_{c_n}$ and $\|u_n-w_n\|_{H^s} \leq
K|1 - \frac{c_n}{c}|$
$$I_{c_n} \leq J(w_n) \leq J(u_n) + L(K)L|1 - \frac{c_n}{c}|,$$
proving that limsup $I_{c_n} \leq \lim J(u_n) = I_c$. This together
with (3.4) imply that
$$\lim_{n\rightarrow \infty} I_{c_n} = I_c .$$
{\bf Lemma 3.2}
\begin{enumerate}
\item If $F$ satisfies (F0) and (F1), then $I_c < 0$ for any $c >
0$.
\item If $F$ satisfies (F1) and (F4), then $I^\infty_c < 0$ for any
$c > 0$.
\end{enumerate}
{\bf Proof} : Let $\varphi$ be a non-negative, radial and radially
decreasing function belonging to $S_c$.\\
Let $0 < \lambda <<< 1$ and set $\varphi_\lambda(x) = \lambda^{N/2}
\varphi(\lambda x)$ then $\varphi_\lambda \in S_c$ and
\begin{eqnarray*}
J(\varphi_\lambda) &=& C_{N,s}
\int\int\;\frac{|\lambda^{N/2}\varphi(\lambda
x)-\lambda^{N/2}\varphi(\lambda y)|^2}{|x-y|^{N+2s}} dxdy\\
&-& \int F(x, \lambda^{N/2} \varphi(\lambda x))dx\\
J(\varphi_\lambda) &\leq& C_{N,s}\int\int \lambda^N
\frac{|\varphi(\lambda x)- \varphi(\lambda y)|^2}{|x-y|^{N+2s}} dx
dy\\
&-& \int_{|x| \geq R} F(x, \lambda^{N/2} \varphi(\lambda x))dx\\
&\leq& C_{N,s} \lambda^{2s} \int\int
\frac{|\varphi(x)-\varphi(y)|^2}{|x-y|^{N+2s}}dx dy\\
&-& \Delta \lambda^{N/2\alpha}\int_{|x| \geq R} |x|^{-p}
\varphi^\alpha(\lambda x)dx\\
&\leq& \lambda^{2s} |\nabla_s \varphi|^2_2 - \Delta
\lambda^{\frac{N}{2}\,\alpha} \lambda^{-N} \lambda^p \int_{|y| \geq
\lambda R} |y|^{-p} \varphi^\alpha(y) dy
\end{eqnarray*}
since $0< \lambda <<< 1$, we certainly have :
\begin{eqnarray*}
J(\varphi_\lambda) &\leq& \lambda^{2s} |\nabla_s\varphi|^2_2 -
\Delta \lambda^{\frac{N}{2}\, \alpha - N+p} \int_{|y| \geq R}
|y|^{-p}\varphi^\alpha(y)dy\\
&\leq& \lambda^{2s} \{C_1 - \lambda^{\frac{N}{2}\, \alpha-N+p-2s}
C_2\}
\end{eqnarray*}
letting $\lambda \rightarrow 0$ and using the fact that $N+2s >
\frac{N}{2} \alpha + p$ the strict  negativity of $I_c$ follows.

b) The proof is dentical.\\
{\bf Lemma 3.3}
\begin{enumerate}
\item If $F$ satisfies (F0), (F1) and (F2) then
$$I_c \leq I_a + I_{c-a}\quad \forall\; a \in (0,c)\eqno{(3.5)}$$
\item If $F$ satisfies (F2), (F4) and (F1) holds true for $F^\infty$
then :
$$I^\infty_c < I^\infty_a + I^\infty_{c-a} \quad \forall\; a \in (0,c)\eqno{(3.6)}$$
\end{enumerate}
{\bf Proof : }
\begin{enumerate}
\item This is a direct consequence of the fact that a real-valied
function $f$ satisfying $f(\theta t) \leq \theta^2 f(t)$ for any
$\theta \geq 1$ does certainly verify :
$$f(c) \leq f(a) + f(a-c)\quad \forall\; a \in (0,c),\;\; [7]$$
\item Following the same steps as in the last part, we can conclude
that : $I^\infty_{\theta c} < \theta^2 I^\infty_c \quad \forall\;
\theta > 1$.\\
Let $c > 0, 0 < a < c$ and $\theta > 1$, we can choose $\varepsilon
> 0$ such that $\varepsilon < - I^\infty_c(1-\theta^{-\sigma})$ and
there exists $v \in S_c$ such that : $I^\infty_c < J^\infty(v) <
I^\infty_c + \varepsilon$ .\\
Hence $$I^\infty_{\theta c} \leq J^\infty(\theta v) \leq
\theta^{\sigma+2} J^\infty(v).$$
 Therefore $I^\infty_{\theta c} \leq
\theta^{\sigma+2} \{I^\infty_c + \varepsilon\}$ $< \theta^{\sigma+2}
I_c^\infty$ by the choice of $\varepsilon$.
\end{enumerate}
{\bf Proof of Theorem 1.2}

Let $(u_n)$ be a minimizing sequence of the problem (1.4).\\
{\bf Vanishing does not occur }:

If it occurs it follows from Lemma I.1 of [7] that
$|u_n|_p\rightarrow 0$ as $n \rightarrow + \infty$ for $p \in (2,
2^\ast_s)$. By (F4)
$$\int F^\infty(x, u_n) \leq B \{|u_n|^{\gamma+2}_{\gamma+2} +
|u_n|^{\ell+2}_{\ell+2}\}.$$ Thus $\displaystyle{\lim_{n\rightarrow
+ \infty}\int} F^\infty(x,u_n) =0$, which implies that liminf
$J^\infty(u_n) \geq 0$,
contradicting the fact that $I^\infty_c < 0$.\\
{\bf Dichotomy does not occur } :

We will use the notation introduced in the appendix :

For $n \geq n_0 : J^\infty(u_n) - J^\infty(v_n) - J^\infty(w_n) =$
$$\frac{1}{2} \int|\nabla_s u_n|^2 - |\nabla_s v_n|^2 -
|\nabla_s w_n|^2 - \int
F^\infty(x,u_n)-F^\infty(x,v_n)-F^\infty(x,w_n)$$
$$\frac{1}{2} \int|\nabla_s u_n|^2 - |\nabla_s v_n|^2 -
|\nabla_s w_n|^2 - \int F^\infty (x,u_n)-F^\infty
(x, v_n + w_n)$$ since supp $v_n \cap$ supp $w_n = \emptyset$
$$\geq - \varepsilon - \int F^\infty(x,u_n) - F^\infty(x, v_n + w_n).$$
Now since $\{v_n\}$ ad $\{w_n\}$ are also bounded in $H^s$, it
follows from the proof of Lemma 3.1 that there exists $C,K > 0$ such
that :
\begin{eqnarray*}
|\int F^\infty(x,u_n) &-& F^\infty(x, v_n + w_n)|\\
&\leq& \sup_{|u|_{H^s \leq K}}|\partial_2
F^\infty(x,u)|_{H^{-s}}|u_n -(v_n + w_n)|_{H^s}\\
&\leq& \sup_{|u|_{H^s}\leq K} |\partial^1_2 F^\infty(x,u)|_{L^2}
|u_n - (v_n+w_n)|_{L^2} \\
&+& \sup_{|u|_{H^s} \leq K} |\partial^2_2 F^\infty(x,u)|_{L^p}
|u_n-(v_n + w_n)|_{L^{p'}}\\
&\leq& C \sup_{|u|_{H^s} \leq K} |u|_{L^2} |u_n-(v_n + w_n)|_{L^2}\\
&+& C \sup_{|u|_{H^s} \leq K} |u|_{L^q}^{1+ \frac{4s}{N}}|u_n - (v_n
+w_n)|_{L^{p'}}\\
&\leq& C_1 K |u_n-(v_n + w_n)|_{L^2} + C_2 K^{1 + \frac{4s}{N}}
|u_n- (v_n + w_n)|_{L^{p'}}
\end{eqnarray*}
so :\\
$J^\infty(u_n) - J^\infty(v_n) - J^\infty(w_n)\geq$\\
$$- \varepsilon - C_1 K|u_n-(v_n+w_n)|_{L^2} + C_2 K^{1 + \frac{4s}{N}}
|u_n-(v_n + w_n)|.$$ Given any $\delta > 0$, we can find
$\varepsilon_\delta \in (0, \delta)$ such that (we have used the
properties of the sequences $(v_n)$ and $(w_n))$
$$J^\infty(u_n)-J^\infty(v_n)-J^\infty(w_n) \geq - \delta.$$
Now let $$a^2_n(\delta) = \int v^2_n$$
$$b^2_n(\delta) = \int w^2_n.$$
Passing  to a subsequence, we may suppose that :
$$a^2_n(\delta)  \rightarrow a^2(\delta)$$
and
$$b^2_n(\delta) \rightarrow b^2(\delta)$$
where
$$|a^2(\delta)-a^2| \leq \varepsilon_\delta < \delta
\mbox{ and } |b^2(\delta)-(c^2-a^2)| < \varepsilon$$ Recalling that
$I^\infty_c$ is continuous, we find that :
\begin{eqnarray*}
I^\infty_c &\geq& \lim_{n\rightarrow \infty} J^\infty(u_n) \geq
\liminf \{J^\infty(v_n) + J^\infty(w_n)\}\\
&\geq& I^\infty_{a(\delta)} + I^\infty_{b(\delta)} - \delta.
\end{eqnarray*}
Letting $\delta$ goes to zero and using again the continuity of
$I^\infty_c$, we obtain :
$$I^\infty_c \geq I^\infty_a + I^\infty_{\sqrt{c^2-a^2}}$$
contradicting Lemma 3.3.\\
Hence {\bf compactness occurs} : so there exists $\{y_n\} \subset
\mathbb{R}^N$ such that for all $\varepsilon > 0$ : $$\int_{B(y_n,
R(\varepsilon))} u^2_n \geq c^2 - \varepsilon.$$ For each $n \in
\mathbb{N}$, we can choose $z_n \in \mathbb{Z}^N$ such that $y_n -
z_n \in [0,1]^N$.\\
Now let $v_n = u_n (x+z_n)$, we certainly have that $|v_n|_{H^s} =
|u_n|_{H^s}$ is bounded and so passing to a subsequence, we may
assume that $(v_n)$ converges weakly to $v$ in $H^s$ in particular
$(v_n)$ converges weakly to $v$ in $L^2$ and  $|v_n|^2_2 = c^2$, but
\begin{eqnarray*}
\int v^2 &\geq& \int_{B(0, R(\varepsilon) + \sqrt{N})} |v|^2\\
&=& \lim_{n\rightarrow \infty} \int_{B(0, R(\varepsilon) +
\sqrt{N})} |v_n|^2 = \lim \int_{B(z_n, R(\varepsilon) + \sqrt{N})}
|v_n|^2
\end{eqnarray*}
and
$$\int_{B(z_n, R\varepsilon)+\sqrt{N})} u^2_n \geq \int_{B(y_n,R(\varepsilon))}
u^2_n \geq c^2 - \varepsilon$$ since $|y_n-z_n| \leq \sqrt{N}$.\\
Hence $|v|^2_{L^2} \geq c^2 - \varepsilon$ $\forall\; \varepsilon >
0$ $\Rightarrow |v|^2_{L^2} \geq c^2$.\\
On the other hand $|v|_2 \leq \liminf |v_n|_2\Rightarrow |v|^2_{L^2}
\leq c^2$.\\
It follows then that $|v|^2_{L^2} = c^2 \Rightarrow |v-v_n|_{L^2}
\rightarrow 0$ as $n \rightarrow \infty$.\\
Furthermore by the periodicity of $F^\infty$ :
$$J^\infty(u_n) = J^\infty(v_n) \rightarrow I^\infty_c$$
and
$$v_n \rightarrow v \mbox{ in } L^p, p \in [2, 2^\ast_s) .$$
It follows that $v_n \rightarrow v$ in $H^s$  and consequently
$\displaystyle{\int}F^\infty(x, v_n) \rightarrow \displaystyle{\int}
F^\infty(x,v)$, which implies that
$$J^\infty(v) = I^\infty_c .$$
{\bf Lemma 3.4}  \\

If $F$ satisfies (F0), (F1), (F2) and (1.4) is achieved then
$$I_c < I_a + I^\infty_{c-a}\quad \forall\; a \in (0,c).$$
 {\bf Proof of Theorem 1.1}

In the following $(u_n)$ is a minimizing sequence of (1.1) and we
will make use of the notation introduced in the appendix.\\
{\bf Vanishing does not occur }:

If it occurs, it would follow from Lemma I.1 of [7] that
$|u_n|_{L^p} \rightarrow 0$ for $p \in (2, 2^\ast_s)$.\\
Combining (F0) and (F3) we have :  For any $\delta > 0, \exists \;
R_\delta > 0$ such that
$$F(x,t) \leq \delta(t^2 + |t|^{\beta+2}) + A'(|t|^{\gamma+2} +|t|^{\ell+2})\quad
\forall\; |x| \geq R_\delta .$$ Hence
$$\int_{|x| \geq R_\delta} F(x, u_n) \leq \delta(|u_n|^2_2 +
|u_n|^{\beta+2}_{\beta+2}) +
A'(|u_n|^{\beta+2}_{\beta+2} + |u_n|^{\ell+2}_{\ell +2}).$$ Thus
$$\limsup_{n\rightarrow \infty}
\int_{|x| \geq R_\delta} F(x, u_n) \leq \delta c^2.$$ On the other
hand :
\begin{eqnarray*}
\int_{|x| \leq R_\delta} F(x,u_n) &\leq& A \int_{|x| \leq R_\delta}
|u_n|^2 + |u_n|^{\ell +2}\\
 &\leq& A \Big\{ |u_n|^{\ell+2}_{\ell+2}
|R_\delta|^{\frac{\ell}{\ell+2}} + |u_n|^{\ell+2}_{\ell+2}\Big\}
\displaystyle{\mathop{\longrightarrow 0}_{n \rightarrow \infty}}
\end{eqnarray*}
Hence for any $\delta > 0$ we have that
$$\limsup_{n\rightarrow \infty} \int F(x,u_n) < \delta c^2$$
and so $\lim \int F(x,u_n) = 0$.\\
But $J(u_n) \rightarrow I_c < 0$ and we obtain the contradiction.\\
{\bf Dichotomy does not occur} :

Suppose that the sequence $\{y_n\}$ is bounded and let us consider :
\begin{eqnarray*}
J(u_n)\!\!\! \!&-&\!\!\!\! J(v_n) - J^\infty(w_n) = \frac{1}{2}
\{|\nabla_s u_n|^2_2 -
|\nabla_s v_n|^2_2 - |\nabla_s w_n|^2_2\}\\
&-& \int F(x,u_n) -  F(x, v_n) - F(x, w_n)\\
&+& \int F^\infty(x, w_n)- F(x,w_n)\\
&\geq& - \varepsilon - \int F(x, u_n) - F(x,v_n + w_n) + \int
F^\infty(x, w_n) - F(x, w_n)
\end{eqnarray*}
since supp $v_n \cap supp w_n = \emptyset$
$$\geq - \varepsilon - \int F(x,u_n) - F(x, v_n + w_n) +
 \int_{|x-y_n| \geq R_n}
F^\infty(x, w_n) - F(x,w_n).$$ Now using the same argument as
before, it follows that :\\
Given $\delta > 0$, we can choose $\varepsilon = \varepsilon_\delta
\in (0, \delta)$ such that $-\varepsilon - \int F(x, u_n) - F(x, v_n
+ w_n) \geq - \delta$ and hence $J(u_n) - J(v_n) - J^\infty(w_n)
\geq - \delta + \int_{|x-y_n| \geq R_n}F^\infty(x, w_n)-F(x,
w_n)$.\\
Given $\eta > 0$, we can find $R > 0$ such that for all $t \in
\mathbb{R}$ and $|x| \geq R$
$$|F^\infty(x,t)-F(x,t)| \leq \eta (t^2 + |t|^{\beta+2}).$$
Now since $R_n \rightarrow \infty$ and we  are supposing that
$\{y_n\}$ is bounded, we have that :
$$\{x : |x-y_n| \geq R_n\} \subset \{x : |x| \geq R\}$$
for $n$ large enough.\\
From this and the boundedness of $w_n$ in $H^s$, it follows that
$$\lim_{n\rightarrow \infty} \int_{|x-y_n| \geq R_n}
 F^\infty(x, w_n)-F(x,w_n) = 0.$$
Now let $$a^2_n(\delta) = \int v^2_n$$
$$b^2_n(\delta) = \int w^2_n.$$
Passing to a subsequence, we may suppose that :
\begin{eqnarray*}
a^2_n (\delta) &\rightarrow& a^2(\delta)\\
b^2_n(\delta) &\mapsto& b^2(\delta)
\end{eqnarray*}
where $|a^2_n(\delta) - a^2| < \varepsilon_\delta < \delta$ and
$|b^2_n(\delta)-(c^2-a^2)|\leq \varepsilon_\delta < \delta$.\\
Recalling that $I_c$ and $I^\infty_c$ are continuous, we find that :
\begin{eqnarray*}
I_c &=& \displaystyle{\lim_{n\rightarrow \infty}}J(u_n) \geq
\liminf\{J(v_n) + J^\infty(w_n)\}-\delta\\
&\geq& \liminf\{I_{a_n(\delta)} + I_{b_n(\delta)}\} - \delta
\end{eqnarray*}
Thus $I_c \geq I_a + I_{\sqrt{c^2-a^2}}-\delta$.\\
Letting $\delta \rightarrow 0$ we get $I_c \geq I_a +
I_{\sqrt{c^2-a^2}}$.\\
Thus the sequence $\{y_n\}$ cannot be bounded and, passing to a
subsequence, we may suppose that $|y_n| \rightarrow \infty$. Now we
obtain a contradiction with Lemma 3.4 by using similar arguments
applied to $J(u_n) -J^\infty(v_n)-J(w_n)$ to show that $I_c \geq
I^\infty_a + I_{\sqrt{c^2-a^2}}$.\\
Thus dichotomy cannot occur and we have {\bf compactness}.\\
According to the appendix, there exists $\{y_n\} \subset
\mathbb{R}^N$ such that $$\int_{B(y_n,R(\varepsilon))} u^2_n \geq
c^2 - \varepsilon \quad \forall\; \varepsilon > 0.$$ Let us first
prove that the sequence $\{y_n\}$ is bounded . If it is not the
case, we may assume that $|y_n| \rightarrow \infty$ by passing to a
subsequence. Now we can choose $z_n \in \mathbb{Z}^N$ such that $y_n
- z_n \in [0,1]^N$.\\
Setting $v_n(x) = u_n(x+z_n)$, we can suppose that $(v_n)$ converges
weakly to $v$ in $H^s$ and  $|v_n-v|_{L^2} \rightarrow 0$ as $n
\rightarrow \infty$ for $2\leq p \leq 2^\ast_s$. Of course
$J^\infty(v_n) = J^\infty(u_n)$.\\
On the other hand $J(u_n)-J^\infty(u_n) =
\displaystyle{\int}F^\infty(x, u_n)-F(x,u_n) = \int F^\infty(x,v_n)
-
F(x-z_n, v_n)$.\\
Now given $\varepsilon > 0$, it follows from (F3) that there exists
$R > 0$ such that :
$$|\int_{x-z_n| \geq R} F^\infty(x,v_n)-F(x-z_n, v_n)| =$$
$$|\int_{|x-z_n| \geq R} F^\infty(x-z_n,v_n) - F(x-z_n, v_n)|\leq$$
\begin{eqnarray*}
&\leq& \varepsilon \int_{|x-z_n|\geq R} |v_n|^2 + |v_n|^{\beta+2}
\leq \varepsilon C\{|v_n|^2_{H^s} + |v_n|_{H^s}^{\beta+2}\\
&\leq& \varepsilon D \mbox{ since } (v_n) \mbox{ is bounded in } H^s
\end{eqnarray*}
On the other hand since $|z_n| \rightarrow \infty$, there exists
$n_R > 0$ such that for all $n \geq n_R$ :
\begin{eqnarray*}
&&|\int_{|x-z_n| \leq R} F^\infty(x,v_n) - F(x-z_n,v_n)| \\
&\leq&\int_{|x| \geq \frac{1}{2} |z_n|} F^\infty(x, v_n) - F(x-z_n,
v_n)| \\
&\leq& A \int_{|x| \geq \frac{1}{2} |z_n|} |v_n|^2 +
|v_n|^{\ell+2}\\
&\leq& A \int_{|x|\geq \frac{1}{2} |z_n|} |v|^2 + |v|^{\ell+2} + A
\int_{|x| \geq \frac{1}{2} |z_n|} |v-v_n|^2 + |v-v_n|^{\ell+2}\\
&\leq& A \int_{|x| \geq \frac{1}{2} |z_n|} |v|^2 + |v|^{\ell+2} + A
\int_{\mathbb{R}^N} [v-v_n|^2 + |v-v_n|^{\ell+2}
\end{eqnarray*}
and hence
$$\lim|\int_{|x-z_n| \geq R_n} F^\infty(x,v_n) - F(x-z_n, v_n)| =0. $$
Thus
$$\liminf \{J(u_n)-J^\infty(u_n)\} \geq - \varepsilon D\quad
\forall\; \varepsilon > 0.$$ And so $I_c = \lim J(u_n) \geq \lim
J^\infty(u_n) \geq I_c^\infty$ contradicting the fact that $I_c <
I^\infty_c$.\\
Hence $\{y_n\}$ is bounded. Set $\rho = \displaystyle{\sup_{n \in
\mathbb{N}}} |y_n|$, it follows that
$$\int_{B(0,R(\varepsilon)+\rho)} u^2_n \geq \int_{B(y_n,R(\varepsilon)}
u^2_n \geq c^2 - \varepsilon \quad \forall\; \varepsilon > 0.$$ Thus
\begin{eqnarray*}
\int u^2 &\geq& \int_{B(0,R(\varepsilon)+\rho)} u^2 =
\lim_{n\rightarrow \infty}\int_{B(0,R(\varepsilon)+\rho)} u^2_n\\
&\geq& c^2 - \varepsilon\quad \forall\; \varepsilon > 0.
\end{eqnarray*}
and hence $\displaystyle{\int}u^2 \geq c^2$, on the other hand $\int
u^2 \leq c^2$. Thus $u \in S_c$ and $|u_n-u|_{L^2} \rightarrow 0$.
By the boundedness of $u_n$ in $H^s$, it follows that
$u_n\rightarrow u$ in $L^p$ for $p\in [2, 2^\ast_s]$, therefore
$$\lim_{n\rightarrow \infty} \int F(x,u_n) = \int F(x,u), \mbox{ implying}$$
that $J(u) = I_c$.
\section*{Appendix} The concentraction compactness Lemma :\\
If $(u_n)$ is a bounded sequence in $H^s$ such that $\int u^2_n =
c^2$, then one of the following alternatives occur.
\begin{enumerate}
\item {\bf Vanishing } : $\displaystyle{\limsup_{y \in \mathbb{R}^N}}
\displaystyle{\int_{y+B_R}}u^2_n = 0$.
\item {\bf Dichotomy} : There exists $a \in (0,c)$ such that $\forall\;
 \varepsilon > 0, \exists\; n_0 \in
\mathbb{N}$ and two bounded sequences in $H^s$ denoted by $v_n$ and
$w_n$ (all depending on $\varepsilon$) such that for every $n \geq
n_0$, we have
$$|\int v^2_n -a^2| < \varepsilon \mbox{ and } |\int w^2_n -(c^2-a^2)| < \varepsilon$$
$$||\nabla_s u_n|^2 - |\nabla_s v_n|^2 - |\nabla_s w_n|^2 \geq -2\varepsilon$$
and
$$|u_n-(v_n-w_n)|_p \leq 4\varepsilon\quad \forall\; p \in [2,2^\ast_s].$$
Furthermore $\exists\; (y_n) \subset \mathbb{R}^N$ and $\{R_n\}
\subset (0,\infty)$ such that $\lim_{n\rightarrow \infty} R_n = +
\infty$ and : $$\left\{ \begin{array}{ll} v_n = u_n &\mbox{ if }
|x-y_n| \leq R_0\\
|v_n| \leq |u_n| &\mbox{ if } R_0 \leq |x-y_n| \leq 2R_0 \\
v_n = 0 &\mbox{ if } |x-y_n| \leq 2R_0
\end{array}\right.$$
$$\left\{ \begin{array}{ll}
w_n = 0&\mbox{ if } |x-y_n| \leq R_n\\
|w_n| \leq |u_n| &\mbox{ if } R_n \leq |x-y_n| \leq 2R_n\\
w_n = v_n &\mbox{ if } |x-y_n| \geq 2R_n
\end{array}\right.$$
with dist (supp$|v_n|$, supp$(w_n))\rightarrow \infty$ as $n
\rightarrow \infty$.
\end{enumerate}
{\bf Compactness } : There exists a sequence $\{y_n\} \subset
\mathbb{R}^N$ such that for all $\varepsilon > 0$, there exists
$R(\varepsilon) > 0$ such that
$$\int_{B(y_n,R(\varepsilon))} u^2_n \geq c^2 - \varepsilon.$$
\section*{References}
\begin{enumerate}
\item C. J. Amick and J. F. Toland, Uniqueness and related analytic
properties for the BenjaminOno equation-a nonlinear Neumann problem
in the plane, Acta Math., 167 (1991), pp. 107-126.
\item L.
Abdelouhab, J. L. Bona, M. Felland, and J.-C. Saut, Nonlocal models
for nonlinear, dispersive waves, Phys. D, 40 (1989), pp. 360-392.
\item  Eleonora Di Nezza G Patalucci, E Valdinocci, Hitchhiker'sguide
To the fractional Sobolev spaces.
\item  M. Fall, personal
communications.
\item H Hajaiej, Variational problems related to
some fractional kinetic equations, Preprint.
 \item H Hajaiej, L
Molinet, T Ozawa, B Wang, Necessary and sufficient conditions for
the fractional Gagliardo-Nirenberg inequalities and applications to
Navier-Stokes and generalized Boson equations, Preprint.
 \item  P.
L. Lions, P.L. Lions, The concentration-compactness principle in the
calculus of variations, the locally compact case, Part 1( p 109-145)
and Part2( p 223-281). Ann Ins H Poincare Vol 1 N4, 1984.
\end{enumerate}
\end{document}